\documentclass[]{siamart171218}      

\usepackage{amsmath}
\usepackage{amssymb}
\usepackage{tikz}
\usepackage{listings}
\newcommand{\R}{{I\!\!R}}
\newcommand{\U}{{\bf u}}

\newcommand{\V}{{\bf v}}
\newcommand{\W}{{\bf w}}
\newcommand{\X}{{\bf x}}

\newcommand{\Y}{{\bf y}}

\newcommand{\ass}{\small \mathrel{{:}{=}}}

\title{Differential Invariants}
\author{Uwe Naumann\thanks{Informatik 12: Software and Tools for Computational Engineering, RWTH Aachen University, Germany. \email{naumann@stce.rwth-aachen.de}}}

\begin{document}

\maketitle


\begin{abstract}
Validation is a major challenge in differentiable programming.
The state of the art is based on algorithmic differentiation. 
Consistency of first-order tangent and adjoint programs is defined by a well-known
first-order differential invariant. This paper generalizes the approach through 
derivation of corresponding differential invariants of arbitrary order.
\end{abstract}

\section{Introduction and State of the Art}

We consider 
implementations of 
multivariate vector functions
\begin{equation} \label{eqn:F}
	F: \R^n \rightarrow \R^m : \X \mapsto \Y=F(\X) 
\end{equation}
over the real (floating-point) numbers $\R$
as sufficiently often continuously differentiable computer programs. 
Distinctly named variables (e.g, $\X$ and $\Y$) are assumed to be {\em unaliased}\footnote{They occupy disjoint system memory locations.} in the given implementation. 
We refer to $F$ as the {\em primal program}. 
Its $\nu$-th derivative (tensor) is denoted as
\begin{equation} \label{eqn:derivs}
	F^{[\nu]}=F^{[\nu]}(\X) \equiv \frac{d^\nu F}{d \X^\nu}(\X) \in \R^{m \times n^\nu} \; , 
\end{equation}
where 
$n^\nu \equiv \overset{\nu~\text{times}}{\overbrace{n \times \ldots \times n}}.$
We set $F'\equiv F^{[1]}$ (Jacobian) and $F''\equiv F^{[2]}$ (Hessian). 
Vectors are printed in bold font. Upper case letters denote other matrices 
and tensors. 
We use $=$ to denote mathematical equality and $\equiv$ in the sense of ``is defined as.''

Index notation $F^{[\nu]}=F^{[\nu]}_{k,j_1,\ldots,j_\nu}$ is used for 
operations on derivative tensors with $k=1,\ldots,m$ and $j_i=1,\ldots,n$ for
$i=1,\ldots,\nu.$ Products with vectors $\U = \U_{j_i}\equiv (u_s) \in \R^n$ are
defined as
\begin{equation} \label{eqn:in1}
F^{[\nu]}_{k,j_1,\ldots,j_\nu} \cdot \U_{j_i} \equiv 
\sum_{s=1}^n F^{[\nu]}_{k,j_1,\ldots,j_{i-1},s,j_{i+1},\ldots,j_\nu} \cdot u_s \; .
\end{equation}
Similarly, 
products with vectors $\W=\W_k\equiv (w_s) \in \R^m$ are defined as
\begin{equation} \label{eqn:in2}
F^{[\nu]}_{k,j_1,\ldots,j_\nu} \cdot \W_k \equiv 
\sum_{s=1}^m F^{[\nu]}_{s,j_1,\ldots,j_\nu} \cdot w_s \; .
\end{equation}
The following rather obvious observation will be used in upcoming proofs.
\begin{lemma} \label{lem1}
For $F$ as in (\ref{eqn:F}) with derivatives as in
(\ref{eqn:derivs})
\begin{enumerate}
\item $F^{[\nu]}_{k,j_1,\ldots,j_\nu} \cdot \U_{j_{i_1}} \cdot \V_{j_{i_2}}
=F^{[\nu]}_{k,j_1,\ldots,j_\nu} \cdot \V_{j_{i_2}} \cdot \U_{j_{i_1}}$
		for $\U,\V \in \R^n,$ and $i_1,i_2 \in \{1,\ldots,\nu\},$ $i_1\neq i_2;$
\item $F^{[\nu]}_{k,j_1,\ldots,j_\nu} \cdot \W_k \cdot \U_{j_i}
=F^{[\nu]}_{k,j_1,\ldots,j_\nu} \cdot \U_{j_i} \cdot \W_k$
for $\W \in \R^m$, $\U \in \R^n$ and $i \in \{1,\ldots,\nu\}.$
\end{enumerate}
\end{lemma}
\begin{proof}
	The lemma follows immediately from the definition of index notation as in (\ref{eqn:in1}) and  (\ref{eqn:in2}) by exploitation of commutativity of scalar addition and multiplication.
\begin{enumerate}
	\item Let $1 \leq i_1<i_2 \leq \nu$ while $\U =\U_{j_{i_1}} \equiv (u_{s_1}) \in \R^n$ and $\V=\V_{j_{i_2}} \equiv (v_{s_2}) \in \R^n$. Then
\begin{align*}
F^{[\nu]}_{k,j_1,\ldots,j_\nu} &\cdot \U_{j_{i_1}} \cdot \V_{j_{i_2}} \\
	&=\sum_{s_2=1}^n \sum_{s_1=1}^n F^{[\nu]}_{k,j_1,\ldots,j_{i_1-1},s_1,j_{i_1+1},\ldots,j_{i_2-1},s_2,j_{i_2+1},\ldots,j_\nu} \cdot u_{s_1} \cdot v_{s_2} \\
&=\sum_{s_1=1}^n \sum_{s_2=1}^n F^{[\nu]}_{k,j_1,\ldots,j_{i_1-1},s_1,j_{i_1+1},\ldots,j_{i_2-1},s_2,j_{i_2+1},\ldots,j_\nu} \cdot v_{s_2} \cdot u_{s_1} \\
&=F^{[\nu]}_{k,j_1,\ldots,j_\nu} \cdot \V_{j_{i_2}} \cdot \U_{j_{i_1}} \; .
\end{align*}
An analogous argument holds for $i_1>i_2.$
\item Similarly, for $\W = \W_k \equiv (w_{s_1}) \in \R^m,$ $\U = \U_{j_i} \equiv (u_{s_2}) \in \R^n$ and $1 \leq i \leq \nu$
\begin{align*}
F^{[\nu]}_{k,j_1,\ldots,j_\nu} &\cdot \W_k \cdot \U_{j_i} \\
&=\sum_{s_2=1}^n \sum_{s_1=1}^m F^{[\nu]}_{s_1,j_1,\ldots,j_{i-1},s_2,j_{i+1},\ldots,j_\nu} \cdot w_{s_1} \cdot u_{s_2} \\
&=\sum_{s_1=1}^m \sum_{s_2=1}^n F^{[\nu]}_{s_1,j_1,\ldots,j_{i-1},s_2,j_{i+1},\ldots,j_\nu} \cdot u_{s_2} \cdot w_{s_1} \\
&=F^{[\nu]}_{k,j_1,\ldots,j_\nu} \cdot \U_{j_i} \cdot \W_k \; .
\end{align*}
\end{enumerate}
\end{proof}
Continuous differentiability up to the required order $\nu$ implies 
partial symmetry of derivative tensors as
$$
F^{[\nu]}_{k,\pi_1(j_1,\ldots,j_\nu)}=
F^{[\nu]}_{k,\pi_2(j_1,\ldots,j_\nu)}
$$
for arbitrary permutations $\pi_1$ and $\pi_2$ of $j_1,\ldots,j_\nu.$

Algorithmic differentiation (AD)
\cite{Griewank2008EDP,Naumann2012TAo}
of the primal program $F$ with respect to (wrt.) $\X$ in 
tangent (also: forward) mode yields the {\em tangent program} 
$$
F^{(1)} : \R^n \times \R^n \rightarrow \R^m : \quad 
(\X,\X^{(1)}) \mapsto \Y^{(1)} = F^{(1)}(\X,\X^{(1)}) \; ,
$$ 
which computes $\Y^{(1)} = \Y^{(1)}_k \in \R^m$ in a given tangent direction $\X^{(1)} = \X^{(1)}_{j_1} \in \R^n$ as
\begin{equation} \label{eqn:T}
		\Y^{(1)}_k \equiv \left [ \frac{d F}{d \X} \right ]_{k,j_1} \cdot \X^{(1)}_{j_1}= F'_{k,j_1} \cdot \X^{(1)}_{j_1} \; .
\end{equation}
Tangents of primal variables are marked by the superscript $\!^{(1)}.$ Tensors are enclosed in square brackets whenever appropriate for clarification of index notation.

AD of the primal program $F$ wrt. $\X$ in 
adjoint (also: reverse) mode yields the {\em adjoint program} 
$$
	F_{(1)} : \R^n \times \R^m \rightarrow \R^n : \quad 
(\X,\Y_{(1)}) \mapsto \X^{(1)} = F_{(1)}(\X,\Y_{(1)}) \; ,
$$
which computes $\X_{(1)} = \X_{(1)_{j_1}} \in \R^n$
in a given adjoint direction $\Y_{(1)} = \Y_{(1)_k} \in \R^m$ as
\begin{equation} \label{eqn:A}
		\X_{(1)_{j_1}} \equiv \left [ \frac{d F}{d \X} \right ]_{k,j_1} \cdot {\Y_{(1)}}_k  = F'_{k,j_1} \cdot {\Y_{(1)}}_k \; .
\end{equation}
Adjoints of primal variables are marked by the subscript $\!_{(1)}.$

AD and adjoint methods in particular play a central role in modern simulation and data science. 
Key applications include computational fluid dynamics \cite{Towara2013TDA}, 
quantitative finance \cite{Giles2006SAF} and 
machine learning \cite{Rummelhart1986Lrb}. A growing number of AD software 
tools have been developed. 
They support a variety of programming languages. 
Coverage includes C/C++ \cite{Griewank1996AAC,SaAlGauTOMS2019}, Fortran \cite{Hascoet2013TTA,Naumann2005AdF} and Matlab \cite{Bischof2002CST,Coleman2016ADi}.
The high level of activity in AD research and development is also documented by so far
seven international conferences on the subject with associated proceedings / special 
post-conference collections; see, for example, \cite{Christianson2018Sio}.
Refer to the AD community's web portal {\tt www.autodiff.org} for further information on research 
groups, software tools and applications. The web presence includes a comprehensive bibliography 
on the subject.

The following well-known first-order differential invariant follows immediately from (\ref{eqn:T}) and (\ref{eqn:A}).
\begin{theorem} \label{the:foi}
For a primal program $F$ as in (\ref{eqn:F}) 
with tangent and adjoint programs evaluating (\ref{eqn:T}) 
and (\ref{eqn:A}) the following first-order differential invariant holds:
	\begin{equation} \label{eqn:foi}
	\X_{(1)_{j_1}} \cdot \X^{(1)}_{j_1}= {\Y_{(1)}}_k \cdot \Y^{(1)}_k \; .
	\end{equation} 
\end{theorem}
\begin{proof}
	$$\X_{(1)_{j_1}} \cdot \X^{(1)}_{j_1}
		=F'_{k,j_1} \cdot {\Y_{(1)}}_k \cdot \X^{(1)}_{j_1} 
		=_{\text{\sc Lemma}~\ref{lem1}} F'_{k,j_1} \cdot \X^{(1)}_{j_1} \cdot {\Y_{(1)}}_k 
		={\Y_{(1)}}_k \cdot \Y^{(1)}_k \;.$$ $\;$
\end{proof}
{\sc Theorem}~\ref{the:foi} can be used to verify the consistency of tangent and 
adjoint programs for given $\X,$ $\X^{(1)}$ and $\Y_{(1)}.$ 
Shared conceptual errors are not detected, for example, 
if the derivative of $\sqrt{x},$ $x>0,$ is incorrectly assumed to be equal to 
$-\frac{1}{2 \cdot \sqrt{x}}.$ (This mistake was made by the author during the implementation of an early prototype of the AD software dco/c++ ... \cite{NAG-TR2-20}.) 
To address this issue the tangent can be approximated
by a finite difference quotient. 
Consistency of finite differences, tangents 
and adjoints increases the likelihood of correctness of the derivative code.
Special care must be taken to control numerical errors inflicted by finite
difference approximation.

\section{Second-Order Differential Invariants} \label{sec:second_order}

Let us generalize the observations from the previous section for second
derivative programs.

\subsection{Second Derivative Programs}

For primal programs as in (\ref{eqn:F}) second derivative programs 
are obtained by differentiation of the first-order 
tangent or adjoint {\em target programs} with respect to $\X$ in tangent or 
adjoint mode.
Four variants of second derivative programs can be generated.

\begin{lemma} \label{lem:TT}
AD of the tangent program
		$$F^{(1)}(\X,\X^{(1)}) \equiv F'_{k,j_1}(\X) \cdot \X^{(1)}_{j_1}= \Y_k^{(1)}$$
wrt. $\X \in \R^n$ in the tangent direction $\X^{(2)} \in \R^n$ yields the {\em tangent of tangent program}
$$
		F^{(1,2)} : \R^n \times \R^n \times \R^n \rightarrow \R^m : \quad 
(\X,\X^{(1)},\X^{(2)}) \mapsto \Y_k^{(1,2)}=F^{(1,2)}(\X,\X^{(1)},\X^{(2)})
$$ 
for computing
\begin{equation}\label{eqn:TT}
		\Y_k^{(1,2)}= F''_{k,j_1,j_2} \cdot \X^{(1)}_{j_1} \cdot \X^{(2)}_{j_2} \; .
\end{equation}
\end{lemma}
Tangents of variables used in the target program (here: the tangent program) are marked by the superscript $\!^{(2)}.$ Chained superscripts are combined as 
$\!^{(1,2)} \equiv {\!^{(1)}}^{(2)}.$
\begin{proof}
\begin{align*}
		\Y_k^{(1,2)}&= \left [ \frac{d F^{(1)}(\X,\X^{(1)})}{d \X} \right ]_{k,j_2} \cdot \X^{(2)}_{j_2} =
		\left [ \frac{d F'_{k,j_1}}{d \X} \right ]_{k,j_1,j_2} \cdot \X^{(1)}_{j_1} \cdot \X^{(2)}_{j_2} \\ &=
F''_{k,j_1,j_2} \cdot \X^{(1)}_{j_1} \cdot \X^{(2)}_{j_2} \; .
\end{align*}
\end{proof}

\begin{lemma} \label{lem:AT}
AD of the tangent program
$$F^{(1)}(\X,\X^{(1)}) \equiv F'_{k,j_1}(\X) \cdot \X^{(1)}_{j_1}= \Y_k^{(1)}$$
wrt. $\X \in \R^n$ in the adjoint direction $\Y^{(1)}_{(2)} \in \R^m$ yields the {\em adjoint of tangent program}
$$
		F^{(1)}_{(2)} : \R^n \times \R^n \times \R^m \rightarrow \R^n : \quad 
(\X,\X^{(1)},\Y^{(1)}_{(2)}) \mapsto \X_{(2)_{j_2}}=F^{(1)}_{(2)}(\X,\X^{(1)},\Y^{(1)}_{(2)})
$$ 
for computing
\begin{equation}\label{eqn:AT}
		\X_{(2)_{j_2}}= F''_{k,j_1,j_2} \cdot \X^{(1)}_{j_1} \cdot \Y^{(1)}_{(2)_k} \; .
\end{equation}
\end{lemma}
Adjoints of variables used in the target (tangent) program carry the subscript $\!_{(2)}.$ 
\begin{proof}
\begin{align*}
		\X_{(2)_{j_2}}&= \left [ \frac{d F^{(1)}(\X,\X^{(1)})}{d \X} \right ]_{k,j_2} \cdot \Y^{(1)}_{(2)_k} =
		\left [ \frac{d F'_{k,j_1}}{d \X} \right ]_{k,j_1,j_2} \cdot \X^{(1)}_{j_1} \cdot \Y^{(1)}_{(2)_k} \\ &=
F''_{k,j_1,j_2} \cdot \X^{(1)}_{j_1} \cdot \Y^{(1)}_{(2)_k} \; .
\end{align*}
\end{proof}

\begin{lemma} \label{lem:TA}
AD of the adjoint program
		$$F_{(1)}(\X,\Y_{(1)}) \equiv F'_{k,j_1}(\X) \cdot \Y_{(1)_k}= \X_{(1)_{j_1}}$$
wrt. $\X \in \R^n$ in the tangent direction $\X^{(2)} \in \R^n$ yields the {\em tangent of adjoint program}
$$
		F^{(2)}_{(1)} : \R^n \times \R^m \times \R^n \rightarrow \R^n : \quad
(\X,\Y_{(1)},\X^{(2)}) \mapsto \X^{(2)}_{(1)_{j_1}}=F^{(2)}_{(1)}(\X,\Y_{(1)},\X^{(2)})
$$ 
for computing
\begin{equation}\label{eqn:TA}
		\X^{(2)}_{(1)_{j_1}}= F''_{k,j_1,j_2} \cdot \Y_{(1)_k} \cdot \X^{(2)}_{j_2} \; .
\end{equation}
\end{lemma}
\begin{proof}
\begin{align*}
		\X^{(2)}_{(1)_{j_1}}&= \left [ \frac{d F_{(1)}(\X,\Y_{(1)})}{d \X} \right ]_{j_1,j_2} \cdot \X^{(2)}_{j_2} =
		\left [ \frac{d F'_{k,j_1}}{d \X} \right ]_{k,j_1,j_2} \cdot \Y_{(1)_k} \cdot \X^{(2)}_{j_2} \\ &=
F''_{k,j_1,j_2} \cdot \Y_{(1)_k} \cdot \X^{(2)}_{j_2} \; .
\end{align*}
\end{proof}

\begin{lemma} \label{lem:AA}
AD of the adjoint program
		$$F_{(1)}(\X,\Y_{(1)}) \equiv F'_{k,j_1}(\X) \cdot \Y_{(1)_k}= \X_{(1)_{j_1}}$$
wrt. $\X \in \R^n$ in the adjoint direction $\X_{(1,2)} \in \R^n$ yields the {\em adjoint of adjoint program}
$$
		F_{(1,2)} : \R^n \times \R^m \times \R^n \rightarrow \R^n : \quad 
(\X,\Y_{(1)},\X_{(1,2)}) \mapsto \X_{(2)_{j_2}}=F_{(1,2)}(\X,\Y_{(1)},\X_{(1,2)})
$$ 
for computing
\begin{equation}\label{eqn:AA}
		\X_{(2)_{j_2}}= F''_{k,j_1,j_2} \cdot \Y_{(1)_k} \cdot \X_{(1,2)_{j_1}} \; .
\end{equation}
\end{lemma}
Adjoints of variables used in the target (adjoint) program are marked by the 
subscript $\!_{(2)}$ as before. Chained subscripts are combined as 
$\!_{(1,2)} \equiv {\!_{(1)_{(2)}}}.$
\begin{proof}
\begin{align*}
		\X_{(2)_{j_2}}&= \left [ \frac{d F_{(1)}(\X,\Y_{(1)})}{d \X} \right ]_{j_1,j_2} \cdot \X_{(1,2)_{j_1}} =
		\left [ \frac{d F'_{k,j_1}}{d \X} \right ]_{k,j_1,j_2} \cdot \Y_{(1)_k} \cdot \X_{(1,2)_{j_1}} \\ &=
F''_{k,j_1,j_2} \cdot \Y_{(1)_k} \cdot \X_{(1,2)_{j_1}} \; .
\end{align*}
\end{proof}

\subsection{Differential Invariants}

\begin{theorem} \label{the:soi1}
For $F$ as in (\ref{eqn:F}), 
$F^{(1,2)}$ as in (\ref{eqn:TT}) and 
$F^{(1)}_{(2)}$ as in (\ref{eqn:AT})
$$
	\X_{(2)_{j_2}} \cdot \X^{(2)}_{j_2} = \Y^{(1,2)}_k \cdot {\Y^{(1)}_{(2)}}_k \; .
$$
\end{theorem}
\begin{proof}
The theorem follows immediately from 
(\ref{eqn:TT}) and (\ref{eqn:AT}) as
	\begin{align*}
{\X_{(2)}}_{j_2} \cdot \X^{(2)}_{j_2} 
&= F''_{k,j_1,j_2} \cdot \X^{(1)}_{j_1} \cdot {\Y_{(2)}^{(1)}}_k \cdot \X^{(2)}_{j_2} \\
&= F''_{k,j_1,j_2} \cdot \X^{(1)}_{j_1}\cdot \X^{(2)}_{j_2} \cdot {\Y_{(2)}^{(1)}}_k 
=\Y^{(1,2)}_k \cdot {\Y^{(1)}_{(2)}}_k \; .
\end{align*}
\end{proof}
{\sc Theorem}~\ref{the:soi1} can be used to verify the consistency of tangent of tangent and adjoint of tangent programs for given $\X,$ $\X^{(1)},$ $\X^{(2)}$ and $\Y^{(1)}_{(2)}.$ 
Tangents can be approximated by finite differences. 
Potentially serious numerical errors should be expected from second-order 
finite difference approximation. Careful tuning of perturbations is 
crucial. High-precision floating-point arithmetic should be considered.

\begin{theorem} \label{the:soi2}
For $F$ as in (\ref{eqn:F}), 
$F^{(1,2)}$ as in (\ref{eqn:TA}) and 
$F^{(1)}_{(2)}$ as in (\ref{eqn:AA})
$$
	\X_{(2)_{j_2}} \cdot \X^{(2)}_{j_2} = \X^{(2)}_{(1)_{j_1}} \cdot {\X_{(1,2)}}_{j_1} \; .
$$
\end{theorem}
\begin{proof}
The theorem follows immediately from 
(\ref{eqn:TA}) and (\ref{eqn:AA}) as
	\begin{align*}
		\X_{(2)_{j_2}} \cdot \X^{(2)}_{j_2} 
&= F''_{k,j_1,j_2} \cdot \Y_{(1)_k} \cdot \X_{(1,2)_{j_1}} \cdot \X^{(2)}_{j_2} \\
&= F''_{k,j_1,j_2} \cdot \Y_{(1)_k} \cdot \X^{(2)}_{j_2} \cdot \X_{(1,2)_{j_1}} 
		= \X^{(2)}_{(1)_{j_1}} \cdot \X_{(1,2)_{j_1}} \; .
\end{align*}
\end{proof}
{\sc Theorem}~\ref{the:soi2} can be used to verify the consistency of tangent of adjoint and adjoint of adjoint programs for given $\X,$ $\Y_{(1)},$ $\X^{(2)}$ and $\X_{(1,2)}.$ 
Tangents can be approximated by finite differences.


\section{Higher-Order Differential Invariants} \label{sec:higher_order}

Derivative programs of order $\nu$ have the form 
$$\V=F^{[\nu]}_{k,j_1,\ldots,j_{\nu}} \cdot V \; ,$$ 
where $\V$ is an indexed tangent or adjoint of $\X$ or $\Y$ and
$V$ denotes a chained (outer) product of 
indexed tangents or adjoints of $\X$ or $\Y.$ For example, $\nu=2$, 
$\V=\X^{(2)}_{(1)_{j_1}}$ and $V=\Y_{(1)_k} \cdot \X^{(2)}_{j_2}$ in 
a tangent of adjoint program. In the following sub- and superscripts of $\V$ 
will be appended to the (possibly empty) chains of sub- and superscripts in
the expression represented by $\V.$ For example,
$$
\V^{(2)}=\begin{cases}
	\Y^{(1,2)}_k & \text{if}~\V=\Y^{(1)}_k \\
	\X^{(2)}_{(1)_{j_1}} & \text{if}~\V=\X_{(1)_{j_1}} 
\end{cases}
\quad
\text{and}
\quad
\V_{(2)}=\begin{cases}
	\Y^{(1)}_{(2)_k} & \text{if}~\V=\Y^{(1)}_k \\
	\X_{(1,2)_{j_1}} & \text{if}~\V=\X_{(1)_{j_1}} 
\end{cases} \; .
$$

\subsection{Higher Derivative Programs}

\begin{theorem} \label{the:hoit}
Let $\Y=F(\X)$ be defined by (\ref{eqn:F}) with first-order tangent and 
adjoint programs defined by (\ref{eqn:T}) and (\ref{eqn:A}). 
AD of a $(\nu-1)$-th derivative program
$$\V=F^{[\nu-1]}_{k,j_1,\ldots,j_{\nu-1}} \cdot V$$ 
in tangent mode yields the $\nu$-th derivative program
$$
\V^{(\nu)}=F^{[\nu]}_{k,j_1,\ldots,j_\nu} \cdot V \cdot \X^{(\nu)}_{j_\nu}
$$
for $\nu \geq 2.$
\end{theorem}
\begin{proof}
	The proof is by induction over the order $\nu$ of differentiation.
	\begin{itemize}
		\item[] $\nu=2:$ 
            Let the first derivative program be the
	 \begin{itemize}
            \item tangent program, 
that is, $\V=\Y^{(1)}_k,$ $F^{[\nu-1]}_{k,j_1,\ldots,j_{\nu-1}}=F^{[1]}_{k,j_1}$ and $V=\X^{(1)}_{j_1}.$  
AD wrt. $\X$ in tangent mode yields 
			 \begin{align*}
				 \V^{(2)}=\Y^{(1,2)}_k&=
\left [ \frac{d (F^{[1]}_{k,j_1} \cdot \X^{(1)}_{j_1})}{d \X}\right ]_{k,j_2} \cdot \X^{(2)}_{j_2} \\
				 &= F^{[2]}_{k,j_1,j_2} \cdot \X^{(1)}_{j_1} \cdot \X^{(2)}_{j_2} =F^{[2]}_{k,j_1,j_2} \cdot V \cdot \X^{(2)}_{j_2}
			 \end{align*}
due to independence of $V=\X^{(1)}_{j_1}$ from $\X.$
            \item adjoint program,
that is, $\V=\X_{(1)_{j_1}}$ $F^{[\nu-1]}_{k,j_1,\ldots,j_{\nu-1}}=F^{[1]}_{k,j_1}$ and
$V=\Y_{(1)_k}.$ AD wrt. $\X$ in tangent mode
yields 
\begin{align*}
	\V^{(2)}=\X_{(1)_{j_1}}^{(2)}&=
\left [ \frac{d (F^{[1]}_{k,j_1} \cdot \Y_{(1)_k})}{d \X}\right ]_{j_1,j_2} \cdot \X^{(2)}_{j_2} \\
			 &=F^{[2]}_{k,j_1,j_2} \cdot \Y_{(1)_k} \cdot \X^{(2)}_{j_2}=F^{[2]}_{k,j_1,j_2} \cdot V \cdot \X^{(2)}_{j_2} 
\end{align*}
due to independence of $V=\Y_{(1)_k}$ from $\X.$
	 \end{itemize}
 \item[] $\nu-1 \Rightarrow \nu:$ AD of
			$\V=F^{[\nu-1]}_{k,j_1,\ldots,j_{\nu-1}} \cdot V$ 
wrt. $\X$ in tangent mode yields
			$$\V^{(\nu)} = \left [ \frac{d ( F^{[\nu-1]}_{k,j_1,\ldots,j_{\nu-1}} \cdot V )}{d \X} \right ]_{k,j_\nu} \cdot \X^{(\nu)}_{j_\nu}
			= F^{[\nu]}_{k,j_1,\ldots,j_{\nu}} \cdot V \cdot \X^{(\nu)}_{j_\nu} 
$$ 
due to independence of $V$ from $\X.$
	\end{itemize}
\end{proof}

\begin{theorem} \label{the:hoia}
Let $\Y=F(\X)$ be defined by (\ref{eqn:F}) with tangent and adjoint programs defined by (\ref{eqn:T}) and (\ref{eqn:A}). AD
of a $(\nu-1)$-th derivative program
$$\V=F^{[\nu-1]}_{k,j_1,\ldots,j_{\nu-1}} \cdot V$$ in adjoint mode 
yields 
the $\nu$-th derivative program
$$
\X_{(\nu)_{j_\nu}}= F^{[\nu]}_{k,j_1,\ldots,j_\nu} \cdot V \cdot \V_{(\nu)} 
$$
for $\nu \geq 2.$
\end{theorem}
\begin{proof}
	The proof is by induction over the order $\nu$ of differentiation.
	\begin{itemize}
		\item[] $\nu=2:$ 
            Let the first derivative program be the
	 \begin{itemize}
            \item tangent program, that is,
$\V=\Y^{(1)}_k,$ $F^{[\nu-1]}_{k,j_1,\ldots,j_{\nu-1}}=F^{[1]}_{k,j_1}$ and
$V=\X^{(1)}_{j_1}.$ AD wrt. $\X$ in adjoint mode yields 
\begin{align*}
	\X_{(2)_{j_2}}&=
	\left [ \frac{d (F^{[1]}_{k,j_1} \cdot \X^{(1)}_{j_1})}{d \X} \right ]_{k,j_2} \cdot \Y^{(1)}_{(2)_k} 
=F^{[2]}_{k,j_1,j_2} \cdot \X^{(1)}_{j_1} \cdot \Y^{(1)}_{(2)_k} \\
	&=F^{[2]}_{k,j_1,j_2} \cdot V \cdot \V_{(2)}
\end{align*}
due to independence of $V=\X^{(1)}_{j_1}$ from $\X.$
            \item adjoint program, that is,
$\V=\X_{(1)_{j_1}}$ $F^{[\nu-1]}_{k,j_1,\ldots,j_{\nu-1}}=F^{[1]}_{k,j_1}$ and
$V=\Y_{(1)_k}.$ AD wrt. $\X$ in adjoint mode
yields 
\begin{align*}
\X_{(2)_{j_2}}&=
	\left [ \frac{d (F^{[1]}_{k,j_1} \cdot \Y_{(1)_k})}{d \X} \right ]_{j_1,j_2} \cdot \X_{(1,2)_{j_1}} 
=F^{[2]}_{k,j_1,j_2} \cdot \Y_{(1)_k} \cdot \X_{(1,2)_{j_1}} \\
	&=F^{[2]}_{k,j_1,j_2} \cdot V \cdot \V_{(2)}
\end{align*}
due to independence of $V=\Y_{(1)_k}$ from $\X.$
	 \end{itemize}
 \item[] $\nu-1 \Rightarrow \nu:$ AD of
			$\V=F^{[\nu-1]}_{k,j_1,\ldots,j_{\nu-1}} \cdot V$ 
wrt. $\X$ in adjoint mode yields
$$
			\X_{(\nu)_{j_\nu}}=
			\left [ \frac{d (F^{[\nu-1]}_{k,j_1,\ldots,j_{\nu-1}} \cdot V)}{d \X}\right ]_{i,j_\nu} \cdot \V_{(\nu)}= F^{[\nu]}_{k,j_1,\ldots,j_{\nu}} \cdot V \cdot \V_{(\nu)}
$$ 
due to independence of $V$ from $\X$ and with the vector $\V$ as the result of the $(\nu-1)$-th 
derivative program being indexed by $i.$
	\end{itemize}
\end{proof}
\paragraph{Examples}
A third derivative program is derived in tangent of adjoint of tangent mode recursively as follows:
\begin{itemize}
	\item[] $\nu=1:$ $\V=\Y_k$, $V=1,$ tangent mode $\Rightarrow$ $\Y^{(1)}_k=F^{[1]}_{k,j_1} \cdot \X^{(1)}_{j_1};$
	\item[] $\nu=2:$ $\V=\Y^{(1)}_k$, $V=\X^{(1)}_{j_1},$ adjoint mode $\Rightarrow$ 
$\X_{(2)_{j_2}}=F^{[2]}_{k,j_1,j_2} \cdot \X^{(1)}_{j_1} \cdot \Y^{(1)}_{(2)_k};$ Note that in this case $i=k$
		with reference to the proof of {\sc Theorem}~\ref{the:hoia}. 
	\item[] $\nu=3:$ $\V=\X_{(2)_{j_2}}$, $V=\X^{(1)}_{j_1} \cdot \Y^{(1)}_{(2)_k}$, tangent mode $\Rightarrow$ 
$$\X_{(2)_{j_2}}^{(3)}=F^{[3]}_{k,j_1,j_2,j3} \cdot \X^{(1)}_{j_1} \cdot \Y^{(1)}_{(2)_k} \cdot \X^{(3)}_{j_3} \; .$$
\end{itemize}
Eight third derivative programs can be generated by application
of tangent mode
\begin{alignat*}{2}
	&\text{to (\ref{eqn:TT})} \Rightarrow\quad&
	\Y_k^{(1,2,3)}&= F^{[3]}_{k,j_1,j_2,j_3} \cdot \X^{(1)}_{j_1} \cdot \X^{(2)}_{j_2} \cdot \X^{(3)}_{j_3} \\
	&\text{to (\ref{eqn:AT})} \Rightarrow&
	\X^{(3)}_{(2)_{j_2}}&= F^{[3]}_{k,j_1,j_2,j_3} \cdot \X^{(1)}_{j_1} \cdot \Y^{(1)}_{(2)_k} \cdot \X^{(3)}_{j_3}  \\
	&\text{to (\ref{eqn:TA})} \Rightarrow&
	\X^{(2,3)}_{(1)_{j_1}}&= F^{[3]}_{k,j_1,j_2,j_3} \cdot \Y_{(1)_k} \cdot \X^{(2)}_{j_2} \cdot \X^{(3)}_{j_3} \\
	&\text{to (\ref{eqn:AA})} \Rightarrow&
	\X^{(3)}_{(2)_{j_2}}&= F^{[3]}_{k,j_1,j_2,j_3} \cdot \Y_{(1)_k} \cdot \X_{(1,2)_{j_1}} \cdot \X^{(3)}_{j_3} 
	\intertext{and of adjoint mode}
	&\text{to (\ref{eqn:TT})} \Rightarrow&
	\X_{(3)_{j_3}}&= F^{[3]}_{k,j_1,j_2,j_3} \cdot \X^{(1)}_{j_1} \cdot \X^{(2)}_{j_2} \cdot \Y_{(3)_k}^{(1,2)} \\
	&\text{to (\ref{eqn:AT})} \Rightarrow&
	\X_{(3)_{j_3}}&= F^{[3]}_{k,j_1,j_2,j_3} \cdot \X^{(1)}_{j_1} \cdot \Y^{(1)}_{(2)_k} \cdot \X_{(2,3)_{j_2}} \\
	&\text{to (\ref{eqn:TA})} \Rightarrow&
	\X_{(3)_{j_3}}&= F^{[3]}_{k,j_1,j_2,j_3} \cdot \Y_{(1)_k} \cdot \X^{(2)}_{j_2} \cdot \X^{(2)}_{(1,3)_{j_1}} \\
	&\text{to (\ref{eqn:AA})} \Rightarrow&
	\X_{(3)_{j_3}}&= F^{[3]}_{k,j_1,j_2,j_3} \cdot \Y_{(1)_k} \cdot \X_{(1,2)_{j_1}} \cdot \X_{(2,3)_{j_2}} \; .
\end{alignat*}

A fifth derivative program is derived in 
		tangent of adjoint of adjoint of tangent of adjoint mode as follows:
\begin{itemize}
	\item[] $\nu=1:$ $\V=\Y_k$, $V=1$, adjoint mode $\Rightarrow$ $\X_{(1)_{j_1}}=F^{[1]}_{k,j_1} \cdot \Y_{(1)_k}$ 
	\item[] $\nu=2:$ $\V=\X_{(1)_{j_1}}$, $V=\Y_{(1)_k}$, tangent mode $\Rightarrow$ $\X_{(1)_{j_1}}^{(2)}=F^{[2]}_{k,j_1,j_2} \cdot \Y_{(1)_k} \cdot \X^{(2)}_{j_2}$
	\item[] $\nu=3:$ $\V=\X_{(1)_{j_1}}^{(2)}$, $V=\Y_{(1)_k} \cdot \X^{(2)}_{j_2}$, adjoint mode $\Rightarrow$ $$\X_{(3)_{j_3}}=F^{[3]}_{k,j_1,j_2,j_3} \cdot \Y_{(1)_k} \cdot \X^{(2)}_{j_2} \cdot \X_{(1,3)_{j_1}}^{(2)}$$
	\item[] $\nu=4:$ $\V=\X_{(3)_{j_3}}$, $V=\Y_{(1)_k} \cdot \X^{(2)}_{j_2} \cdot \X_{(1,3)_{j_1}}^{(2)}$, adjoint mode $\Rightarrow$ 
$$\X_{(4)_{j_4}}=F^{[4]}_{k,j_1,j_2,j_3,j_4} \cdot \Y_{(1)_k} \cdot \X^{(2)}_{j_2} \cdot \X_{(1,3)_{j_1}}^{(2)} \cdot \X_{(3,4)_{j_3}}$$
	\item[] $\nu=5:$ $\V=\X_{(4)_{j_4}}$, $V=\Y_{(1)_k} \cdot \X^{(2)}_{j_2} \cdot \X_{(1,3)_{j_1}}^{(2)} \cdot \X_{(3,4)_{j_3}}$, tangent mode $\Rightarrow$ 
$$\X_{(4)_{j_4}}^{(5)}=F^{[5]}_{k,j_1,j_2,j_3,j_4,j_5} \cdot \Y_{(1)_k} \cdot \X^{(2)}_{j_2} \cdot \X_{(1,3)_{j_1}}^{(2)} \cdot \X_{(3,4)_{j_3}} \cdot \X^{(5)}_{j_5}.$$
\end{itemize}

\subsection{Differential Invariants}

\begin{theorem}
	For $\nu$-th derivative programs as in {\sc Theorems}~\ref{the:hoit} and \ref{the:hoia} 
	$$
	\X_{(\nu)_{j_\nu}} \cdot \X^{(\nu)}_{j_\nu}=
	\V_{(\nu)} \cdot \V^{(\nu)} \; .
	$$
\end{theorem}
\begin{proof}
	Let $ \V=F^{[\nu-1]}_{k,j_1,\ldots,j_\nu-1} \cdot V.$
With
	$$\X_{(\nu)_{j_\nu}}=F^{[\nu]}_{k,j_1,\ldots,j_\nu} \cdot V \cdot \V_{(\nu)} \;\;
	\text{and} \;\; \V^{(\nu)}=F^{[\nu]}_{k,j_1,\ldots,j_\nu} \cdot V \cdot \X^{(\nu)}_{j_\nu}$$
it follows that 
	\begin{align*}
		\X_{(\nu)_{j_\nu}} \cdot \X^{(\nu)}_{j_\nu}&=
			F^{[\nu]}_{k,j_1,\ldots,j_\nu} \cdot V \cdot \V_{(\nu)} \cdot \X^{(\nu)}_{j_\nu} \\
			&= F^{[\nu]}_{k,j_1,\ldots,j_\nu} \cdot V \cdot \X^{(\nu)}_{j_\nu} \cdot \V_{(\nu)} 
			= \V^{(\nu)} \cdot \V_{(\nu)} 
			= \V_{(\nu)} \cdot \V^{(\nu)} 
	\end{align*}
\end{proof}
Note that the $2^{\nu-1}$ derivative programs of order $\nu-1$ yield $\nu$ 
distinct $\V.$ Hence, 
$\nu$ differential invariants can be derived. For example, 
\begin{itemize}
	\item[] $\V=\Y^{(1,2)}_k \Rightarrow \X_{(3)_{j_3}} \cdot \X^{(3)}_{j_3} = \Y^{(1,2)}_{(3)_k} \cdot \Y^{(1,2,3)}_k$
	\item[] $\V=\X_{(2)_{j_2}} \Rightarrow \X_{(3)_{j_3}} \cdot \X^{(3)}_{j_3} = \X_{(2,3)_{j_2}} \cdot \X^{(3)}_{(2)_{j_2}}$
	\item[] $\V=\X^{(2)}_{(1)_{j_1}} \Rightarrow \X_{(3)_{j_3}} \cdot \X^{(3)}_{j_3} = \X^{(2)}_{(1)_{j_1}} \cdot \X^{(2,3)}_{(1)_{j_1}}$
\end{itemize}
for $\nu=3.$ One out of six sixth-order differential invariants over the
$2^6=64$ different sixth derivative programs is
$
	\X_{(6)_{j_6}} \cdot \X^{(6)}_{j_6} = \X^{(5)}_{(4,6)_{j_4}} \cdot \X^{(5,6)}_{(4)_{j_4}} .
$

\section{Discussion}

Any AD tool capable of generating derivative programs of arbitrary order can be used to
implement the validation of differential invariants. Source code transformation tools such as
Tapenade \cite{Hascoet2013TTA} need to be applicable to their own output. While this is typically straight
forward for tangent of $\ldots$ of tangent programs the repeated application of adjoint mode may
cause difficulties due to technical details specific to the given AD tool. 

Overloading tools such as dco/c++ \cite{NAG-TR2-20} need to allow for recursive instantiation of their
derivative types with derivative types of lower order; dco/c++ supports this feature through
nested C++ templates. 
Derivative programs of arbitrary order can be generated by arbitrary nesting of tangent and 
adjoint types.
From a practical perspective this level of flexibility may not be crucial. Higher-order adjoint
programs can always be generated as tangent of $\ldots$ of tangent of adjoint programs provided
continuous differentiability of the primal program up to the required order.

Differential invariants can be used as a debugging criterion for derivative code generated by AD.
The primal program $F$ evaluates a partially ordered sequence of
differentiable {\em elemental functions}
$\Phi_s=\Phi_s(\V_r)_{r \prec s} : \R^{n_s} \rightarrow \R^{m_s}$
as a {\em single assignment code}\footnote{Each variable is assumed to be written once.}
$$
\V_s \ass \Phi_s(\V_r)_{r \prec s} \quad
\text{for}~s \ass 1,\ldots,q
$$
and where, adopting the notation from
\cite{Griewank2008EDP}, $r \prec s$ if and only if $\V_r$ is an argument of $\Phi_s.$ 
We use $\ass$ to denote assignment as defined by imperative programming languages.

AD of $F$ results in the augmentation of the latter with code for computing tangents or/and adjoints.
For example, AD of the single assignment code in tangent mode yields the tangent single 
assignment code
$$
	\left .
	\begin{split}
		\W \equiv \V_s&\ass \Phi_s(\V_r)_{r \prec s} \equiv \Phi_s(\U) \\
		\W^{(1)}_k&\ass \left [ \Phi'_s \right ]_{k,j_1} \cdot \U^{(1)}_{j_1} \\
	\end{split}
	\quad \right \rbrace \quad \text{for}~s=1,\ldots,q \; .
	$$
By the chain rule of differentiation the resulting tangent program computes
\begin{equation} \label{eqn:tangent}
		\Y= F(\X); \quad
		\Y^{(1)}_k = F'_{k,j_1} \cdot \X^{(1)}_{j_1} \; .
\end{equation}
Given values for the inputs $\X$ and $\X^{(1)}$ yield values for both outputs $\Y$ and $\Y^{(1)}.$
Obviously, (\ref{eqn:T}) is contained within  (\ref{eqn:tangent}). The adjoint single assignment code
can be derived analogously.

Stepping forward through a single assignment code with support for the propagation of
both tangents and adjoints enables debugging of derivative code as follows: Initialization of
$\X^{(1)}$ (for example randomly) in addition to $\X$ yields $\V^{(1)}_s$ for $s=1,\ldots,q.$
For each (or selected) $s$ the initialization of $\V_{s_{(1)}}$ followed by backward propagation
of adjoints yields $\X_{(1)}.$ Consistency of tangents and adjoints up to the current $s$ can
be validated by checking the differential invariant $\X_{(1)_{j_1}} \cdot \X^{(1)}_{j_1} = [\V_{s_{(1)}}]_k \cdot [\V^{(1)}_s]_k.$ Additional evidence for the desirable correctness of the adjoint 
program can be obtained by approximating the tangents by finite differences.

Most established AD software tools support user-defined elemental functions. Their built-in
elemental functions can typically be expected to be correct leaving user intervention as
the most likely source of errors. The sketched debugging algorithm 
enables the localization and subsequent correction of potential errors.

The formalism extends seamlessly to higher derivative programs. Implementation in the context of
AD software yields a number of technical challenges the discussion of which is beyond the scope 
of this paper.

\section{Conclusion}

AD as a form of differentiable programming has become an indispensable ingredient of 
state of the art numerical methods. Software tools for AD provide valuable support for
the (semi-)automatic generation of derivative programs. Validation of correctness and debugging
of such programs poses a serious challenge. The work presented in this paper aims to set the 
mathematical stage for the development
of corresponding methods and for their highly desirable implementation.


\end{document}